\documentclass{amsart}

\usepackage{amsmath,amssymb,amsfonts,amstext,amsthm}
\usepackage{url}
\usepackage{graphicx}
\input xy
\xyoption{all}
\usepackage[all]{xy}

\usepackage{subfiles}
\usepackage{tikz}

\newcommand{\nwc}{\newcommand}

\nwc{\aaa}{\mathcal{F}}
\nwc{\aap}{\mathcal{F}_{P}}
\nwc{\al}{\alpha}

\nwc{\C}{\mathbb{C}}
\nwc{\cb}{\overline{C}}
\nwc{\ccc}{\mathfrak{c}}
\nwc{\ch}{\widehat{C}}
\nwc{\cin}{\textbf{(v)}}
\nwc{\cl}{C'}
\nwc{\cp}{\mathcal{C}_{P}}
\nwc{\cpll}{\mathfrak{c}_{P'}}
\nwc{\ct}{\widetilde{C}}

\nwc{\dd}{\mathcal{L}}
\nwc{\ddd}{\mathfrak{d}}
\nwc{\ddl}{\mathcal{L}'}
\nwc{\dlp}{\delta_{P}}
\nwc{\doi}{\textbf{(ii)}}

\nwc{\enq}{$$}

\nwc{\fl}{\flushleft}
\nwc{\fff}{\mathcal{F}}
\nwc{\ffp}{\mathcal{F}_{P}}
\nwc{\ffq}{\mathcal{F}_{Q}}
\nwc{\ffl}{\mathcal{F}'}

\nwc{\G}{\mathcal{G}}
\nwc{\Ga}{\Gamma}
\nwc{\gtl}{\widetilde{g}}

\nwc{\hra}{\hookrightarrow}
\nwc{\hua}{h^{1}(C,\aaa )}

\nwc{\kk}{{\rm K}}

\nwc{\llb}{\mathcal{L}}

\nwc{\mb}{\mathbb}
\nwc{\mc}{\mathcal}
\nwc{\mm}{\mathfrak{m}}
\nwc{\mmp}{\mathfrak{m}_{P}}
\nwc{\mpd}{\mathfrak{m}_{P}^{2}}

\nwc{\nn}{\mathbb{N}}

\nwc{\ob}{\overline{\mathcal{O}}}
\nwc{\obr}{\mathcal{O}^*}
\nwc{\obp}{\overline{\mathcal{O}}_P}
\nwc{\och}{\mathcal{O}_{\hat{C}}}
\nwc{\oh}{\hat{\mathcal{O}}}
\nwc{\ohp}{\hat{\mathcal{O}}_{P}}
\nwc{\ol}{\mathcal{O}'}
\nwc{\oma}{\Omega (\mathfrak{a})}
\nwc{\omo}{\Omega (\mathcal{O})}
\nwc{\oo}{\mathcal{O}}
\nwc{\op}{\mathcal{O}_P}
\nwc{\opc}{\mathcal{O}_{P,C}}
\nwc{\oph}{\hat{\mathcal{O}}_{P}}
\nwc{\opl}{\mathcal{O}_{P}'}
\nwc{\oplc}{\mathcal{O}_{P,C}'}
\nwc{\opll}{\mathcal{O}_{P'}}
\nwc{\opt}{\tilde{\mathcal{O}}_{P}}
\nwc{\optt}{{\mathcal{O}}_{\tilde{P}}}
\nwc{\oq}{\mathcal{O}_{Q}}
\nwc{\oqt}{\tilde{\mathcal{O}}_{Q}}
\nwc{\ot}{\widetilde{\mathcal{O}}}
\nwc{\overop}{\bar{\oo}_{P}}

\nwc{\pb}{\overline{P}}
\nwc{\pbb}{P^*}
\nwc{\pbi}{\overline{P_{i}}}
\nwc{\pbr}{\overline{P_{r}}}
\nwc{\pgmd}{\mathbb{P}^{g+2}}
\nwc{\pgmu}{\mathbb{P}^{g+1}}
\nwc{\ph}{\hat{P}}
\nwc{\pp}{\mathbb{P}}
\nwc{\prv}{\noindent\textbook{Proof}:}
\nwc{\pt}{\widetilde{P}}
\nwc{\ptl}{\tilde{P}}
\nwc{\pum}{\mathbb{P}^{1}}

\nwc{\qh}{\hat{Q}}
\nwc{\qtl}{\tilde{Q}}
\nwc{\qua}{\textbf{(iv)}}

\nwc{\ra}{\rightarrow}
\nwc{\rh}{\hat{R}}

\nwc{\sei}{\textbf{(vi)}}
\nwc{\sep}{\beq\ast\ \ast\ \ast\enq}
\nwc{\sig}{\sigma}
\nwc{\Sig}{\Sigma}
\nwc{\ssp}{S_{P}}
\nwc{\sss}{{\rm S}}

\nwc{\tre}{\textbf{(iii)}}

\nwc{\um}{\textbf{(i)}}

\nwc{\vpb}{v_{\overline{P}}}
\nwc{\vtxp}{\widetilde{V}_{x,P}}
\nwc{\vxp}{V_{x,P}}

\let \wt=\widetilde
\let \La=\Lambda
\nwc{\wh}{\hat{\omega}}
\nwc{\whp}{\hat{\omega}_{P}}
\nwc{\woch}{\omega\cdot\mathcal{O}_{\hat{C}}}
\nwc{\woh}{\omega\cdot\hat{\mathcal{O}}}
\nwc{\ww}{\omega}
\nwc{\wwb}{\omega^*}
\nwc{\wwct}{\omega _{\widetilde{C}}}
\nwc{\wwh}{\widehat{\omega}}
\nwc{\wwhp}{\widehat{\omega}_P}
\nwc{\wwp}{\omega _{P}}
\nwc{\wwt}{\widetilde{\omega}}
\nwc{\wwtp}{\widetilde{\omega}_P}

\nwc{\zz}{\mathbb{Z}}

\newtheorem{coro}{Corollary}[section]

\newtheorem{dfn}[coro]{Definition}

\newtheorem{prop}[coro]{Proposition}

\newtheorem{rem}[coro]{Remark}
\newtheorem{rems}[coro]{Remarks}
\newtheorem{thm}[coro]{Theorem}

\let \fl=\flushleft

\let \ga=\gamma
\let \sub=\subset
\let \be=\beta
\let \al=\alpha

\begin{document}

\title{Dimension counts for cuspidal rational curves via semigroups}

\author{Ethan Cotterill}
\address{Instituto de Matem\'atica, UFF
Rua M\'ario Santos Braga, S/N,
24020-140 Niter\'oi RJ, Brazil}
\email{cotterill.ethan@gmail.com}

\author{Lia Feital}
\address{Departamento de Matem\'atica, CCE, UFV
Av. P H Rolfs s/n, 36570-000 Vi\c{c}osa MG, Brazil}
\email{liafeital@ufv.br}

\author{Renato Vidal Martins}
\address{Departamento de Matem\'atica, ICEx, UFMG
Av. Ant\^onio Carlos 6627,
30123-970 Belo Horizonte MG, Brazil}
\email{renato@mat.ufmg.br}




\maketitle
\vspace{-20pt}
\begin{abstract}
We study singular rational curves in projective space, deducing conditions on their parameterizations from the value semigroups $\sss$ of their singularities. We prove that a natural heuristic based on nodal curves for the codimension of the space of nondegenerate rational curves of arithmetic genus $g>0$ and degree $d$ in $\mb{P}^n$, viewed as a subspace of all degree-$d$ rational curves in $\mb{P}^n$, holds whenever $g$ is small. On the other hand, we show that this heuristic fails in general, by exhibiting an infinite family of examples of Severi-type varieties of rational curves containing ``excess" components of dimension strictly larger than the space of $g$-nodal rational curves.
\end{abstract}

\section*{Introduction}

Rational curves are essential tools for classifying complex algebraic varieties. Establishing dimension bounds for families of embedded rational curves that admit (analytic isomorphism classes of) singularities of a particular type arises naturally in this context; see, for example \cite{JK} and \cite{Co}, where such bounds are used to infer (dimension-theoretic) information about parameter spaces of rational curves embedded in general hypersurfaces. It is also a basic fact \cite{Lau} that any curve singularity occurs along a rational curve.

\medskip
For our purposes, it will be most useful to view a rational curve as a morphism: given a choice of ambient dimension $n$ and degree $d \geq n$, any such is the image of an $(n+1)$-tuple of holomorphic functions $f=(f_0,\dots, f_n)$ with $f_i \in H^0(\mc{O}_{\mb{P}^1}(d))$. Explicitly, we may write
\begin{equation}\label{parametrizacao}
f_i= a_{i,0} t^d+ a_{i,1} t^{d-1}u+ \cdots+ a_{i,d}u^d
\end{equation}
where $t, u$ are homogeneous coordinates on $\mb{P}^1$, and $a_{i,j}$ are complex coefficients for every $i=0, \dots,n$ and $j=0,\dots,d$. From this point of view, the natural parameter space for rational curves is thus the Grassmannian $\mb{G}(n,d)$.

\medskip
Viewing rational curves as points of the Grassmannian is the point of departure for Griffiths and 
Harris' proof \cite{GH} of the {\it Brill--Noether theorem} for general curves of genus $g \geq 2$. In that work, they show that rational curves with $g$ cusps serve as dimension-theoretic surrogates for general curves of genus $g$, in the sense that the dimensions of the spaces of linear series on such curves are the (generically) expected ones. Eisenbud and Harris subsequently reworked (and simplified the proof of) the Brill--Noether result using rational curves with $g$ cusps \cite{EH}. The fact that cuspidal, as opposed to nodal, rational curves are the most-amenable to dimension estimates in the Grassmannian is a recurring theme in this work as well: the {\it value semigroups} of unibranch, i.e. cuspidal, singularities have a much simple structure than those of singularities with an arbitrary number of branches.

\medskip
A key invariant of cusps is their {\it ramification}, or inflection, in a point along their underlying curve. In each point $P \in \mb{P}^1$, the morphism $f: \mb{P}^1 \ra \mb{P}^n$ is determined by local sections $(\sig_0,\dots,\sig_n)$ vanishing at orders $a_0 < \dots < a_n$ in $P$. Generically, we have $(a_0,a_1,\dots,a_n)= (0,1,\dots,n)$, and the deviation 
\begin{equation}\label{ramif_profile}
{\bf \al}= (a_0,a_1,\dots,a_n)-(0,1,\dots,n)
\end{equation}
indexes a Schubert variety $\Sig_{\bf \al}\subset\mb{G}(n,d)$ of codimension $|\al|$ associated to certain incidence conditions with respect to an osculating flag of a rational normal curve in $\mb{P}^d$ (and canonically specified by \eqref{ramif_profile}). Each such Schubert variety describes those centers of projection that map rational normal curves in $\mb{P}^d$ to unicuspidal degree-$d$ rational curves in $\mb{P}^n$; for an interesting recent extension to the setting of rational curves with multibranch singularities, see \cite{BIV}. A crucial fact is that $m$ Schubert varieties $\Sig_{\bf \al_i}$, for $1 \leq i \leq m$, obtained by ramification in $m$ distinct points $P_i \in \mb{P}^1$ intersect dimensionally transversely; i.e., their codimensions are additive. 

\medskip
The ramification of a cusp imposes linear conditions on the coefficients of the underlying morphism $f=(f_0,\dots,f_n): \mb{P}^1 \ra \mb{P}^n$, namely the vanishing of partial derivatives of the $f_i, i=0,\dots, n$. In general, a cusp is not specified by its ramification. However, the semigroup of a cusp provides a natural set of additional, non-linear conditions on coefficients of the $f_i$, which arise from the multiplicative structure of the local ring of the singularity. At the outset of this project, it seemed reasonable to conjecture that, taken together, these two sets always impose at least $(n-2)g$ independent conditions on morphisms; this was the position we took in \cite{CFM}. Indeed, $(n-2)g$ is the codimension inside the space $M^n_d$ of morphisms $\mb{P}^1 \ra \mb{P}^n$ of degree $d \gg g$ comprising maps with $g$-nodal images; and Joe Harris' celebrated proof of the irreducibility of the Severi variety \cite{H1} shows that all degree-$d$ plane curves of arithmetic genus $g$ lie in the closure of the locus of $g$-nodal rational plane curves.

\medskip
It turns out, however, that when the ambient dimension $n \geq 3$ is sufficiently large, there exist rational unicuspidal curves of arithmetic genus $g$ that lie outside of the closure of the locus of $g$-nodal rational curves. Indeed, we construct infinite families of examples of such curves explicitly; each of our examples is associated with a parameterization of a cusp with generic coefficients, whose ramification determines an arithmetic sequence of a particular type. The upshot is that the corresponding Severi-type varieties $M^n_{d,g}$ are reducible, in the same way that Hilbert schemes of finite-length subschemes of a fixed ambient space of dimension at least 3 are typically reducible.

\medskip
Our method for producing dimension estimates for singular rational curves, which uses the natural stratification of singularities by semigroups, should be viewed as an important first step in the direction of systematically classifying these objects. In subsequent work, we will expand upon this classification scheme in two ways. First, we will show \cite{CLM} how a more systematic implementation of our method naturally yields explicit estimates (which in many cases are sharp) for the codimension of the subvarieties of $M^n_d$ associated with unicuspidal rational curves whose associated value semigroups are {\it $\gamma$-hyperelliptic} in the sense of \cite{To2}. Second, we will show \cite{CFHM} how to systematically implement our method in order to count conditions imposed by {\it multibranch} singularities. The task of organizing and systematically enumerating value semigroups associated with multibranch singularities, as there is no obvious native structure that indexes all of these.

\subsection{Roadmap}
Our first result, Thm~\ref{unibranch_theorem}, establishes that the expected codimension $(n-2)g$ is achieved for morphisms whose images have at-worst cusps, provided that $g \leq 8$. The proof uses the stratification of singularities according to to their value semigroups; in the unibranch case, these are {\it numerical semigroups}, i.e. subsemigroups of $\mb{N}$, and thus form a tree. Theorem~\ref{unibranch_theorem} generalizes some of the results obtained in the papers \cite{JK}, \cite{Co1}, and \cite{Co}, where they are applied towards dimension-counting for rational curves on very general hypersurfaces. A priori there are no obstructions to extending the validity of our results to singularities of higher genus; however, the number of semigroups to be analyzed grows rapidly with $g$.

\medskip
On the other hand, our second result, Thm~\ref{severi_counterexamples}, shows that in general the expected codimension $(n-2)g$ fails to hold. We exhibit an infinite number of examples of this phenomenon associated to unicuspidal rational curves whose value semigroups are of a particular {\it $\ga^{\ast}$-hyperelliptic} type.

\medskip

\subsection{Conventions}
We work over $\mb{C}$. By {\it rational curve} we always mean a projective curve of geometric genus zero. We denote by $M^n_d$ the space of nondegenerate morphisms $f: \mb{P}^1 \ra \mb{P}^n$ of degree $d>0$. Here each morphism is identified with the set of coefficients of its homogeneous parametrizing polynomials, so $M^n_d$ is a space of frames over $\mb{G}(n,d)$. We denote by $M^n_{d,g} \subset  M^n_d$ the subvariety of morphisms whose images have arithmetic genus $g>0$. These curves are necessarily singular. Clearly, $M^n_{d,g}$ contains all curves with $g$ simple nodes or $g$ simple cusps.

\medskip
It will also be useful to consider parameter spaces associated to a fixed choice of value semigroup ${\rm S}$. Accordingly, let $\mc{V}_{\rm S} \sub M^n_{d,g}$ denote the space of nondegenerate morphisms $f: \mb{P}^1 \ra \mb{P}^n$ with (fixed) degree $d>0$, image of arithmetic genus (at least) $g$, and a singularity with value semigroup ${\rm S}$ of genus (exactly) $g$. Here the {\it genus} of a value semigroup encodes the contribution of the underlying singularity to the arithmetic genus of the underlying projective curve. We use $a_{i,j}$ to denote the coefficient $[t^iu^{d-i}]f_j$ of $t^i u^{d-i}$ in $f_j$, $j=0, \dots, n$. We reserve the letter $c$ for the {\it conductor} of a semigroup ${\rm S}$. 

\subsection{A heuristic for dimension counts}
Requiring a morphism $f: \mb{P}^1 \ra \mb{P}^n$ to map distinct points $p_1, p_2 \in \mb{P}^1$
to the same image $q \in \mb{P}^n$ imposes $2n$ linear conditions on the coefficients of $f$.
Allowing the preimages and image to vary yields $(n-2)$ linear conditions. Since a simple node
has arithmetic genus 1, it might seem reasonable to expect more generally that singularities of
arithmetic genus $g$ impose at least $(n-2)g$ conditions on morphisms $f: \mb{P}^1 \ra \mb{P}^n$.
In other words, we'd expect that ${\rm cod}(M^n_{d,g},M^n_d)\geq (n-2)g$, when $d$ is
sufficiently large relative to $g$.

\medskip
Note that $(n-2)g$ is the {\it actual} codimension of the is the locus of $g$-nodal rational curves of degree $d \gg g$ in $\mb{P}^n$. Indeed, the coefficient matrix associated with any $g$ fixed choices of nodal targets together with their $2g$ preimages in a Vandermonde matrix, whose determinant is nonvanishing. 

\subsection*{Acknowledgements} We would like to thank Maksym Fedorchuk, Joe Harris, Nathan Ilten, Steve Kleiman, Karl-Otto St\"ohr, Fernando Torres, and Filippo Viviani for illuminating conversations, as well as the mathematics department at UFMG for making this collaboration possible. Special thanks are due to Nicola Maugeri and Giuseppe Zito, who produced the first counterexample to what previously was our ``$(n-2)g$ conjecture".The first and third authors were partially supported by CNPq grant numbers 309211/2015-8 and 306914/2015-8, respectively. The second author is supported by FAPEMIG.

\section{Counting conditions for unibranch singularities}
Unibranch singularities form a naturally distinguished (simple) class of singularities. Accordingly, it makes sense to ask for dimension estimates for rational curves with at-worst-unibranch singularities. We will prove that our na\"ive codimension estimate holds for rational curves of arithmetic genus at most 8:

\begin{thm}\label{unibranch_theorem}
Let $\mathcal{V}:= \bigcup_{{\rm S} \sub \mb{N}} \mc{V}_{\rm S} \subset M^n_{d,g}$ be the subvariety consisting of rational curves with at-worst-unibranch singularities. 
Suppose, moreover, that $g \leq 8$ and $d \geq \max(n,2g-2)$. 
Then
\[
{\rm cod}(\mathcal{V},M^n_d)\geq (n-2)g.
\]
\end{thm}

\begin{rems} 

The choice of genus threshold in Thm~\ref{unibranch_theorem} is essentially arbitrary (indeed, our codimension estimate seems to hold in every computable case) but our method of proof requires us to consider a certain number of exceptional cases, whose number grows as the genus increases.

On the other hand, the condition $d \geq n$ is imposed by the requirement that our rational curves be nondegenerate. It is less clear what a reasonable lower threshold for the degree as a function of the genus should be, but the assumption that $d \geq 2g-2$ is well-adapted to the analysis of conditions beyond ramification in the proof to follow; it also includes the (canonical) case in which $d=2g-2$ and $n=g-1$.

\end{rems}

\medskip
The proof of Thm~\ref{unibranch_theorem} invokes a number of standard tools from linear series and singularities, which we review now. Accordingly, let $P\in C:=f(\pum)\subset\mathbb{P}^n$ be a unibranch singularity. Then $P$ admits a local parametrization $\psi: t \mapsto (\psi_1(t), \dots, \psi_n(t))$ corresponding to a map of rings
\begin{gather*}
\begin{matrix}
\phi : & R:= \mb{C}[[x_1,\dots,x_n]] & \longrightarrow & \mb{C}[[t]]\\
         &x_i            & \longmapsto     &  \psi_i(t)
\end{matrix}
\end{gather*}
Let $v:\mb{C}[[t]] \ra \mb{N}$ denote the standard valuation induced by the assignment $t \mapsto 1$. Let ${\rm S}:=v(\phi(R))$ denote the numerical {\it value semigroup} of $P$. The (local) {\it genus} of the singularity at $P$ is $\delta_P:=\#(\mb{N}\setminus{\rm S})$, and the (global arithmetic) genus of $C$ is the sum of all of these local contributions:
$$
g=\sum_{P\in C}\delta_P.
$$
It will be convenient in what follows to think of the genus of a unibranch singularity as an invariant of the associated numerical semigroup ${\rm S}$; the definition remains unchanged. 
The {\it weight} of a unibranch singularity (or equivalently, of its semigroup ${\rm S}$) of genus $g$ is $W_{\rm S}= \sum_{\ell \in \mb{N}\setminus {\rm S}} \ell - \frac{g(g+1)}{2}$. Here $G_{\rm S}:= \mb{N}\setminus {\rm S}$ is the {\it gap set} of ${\rm S}$.

\subsection{Beginning of the proof of Thm~\ref{unibranch_theorem}}
Let $m_i$ denote the $i$th positive integer in ${\rm S}$ (ordered from smallest to largest), and let $a_0 < a_1 <\cdots < a_n$ be the vanishing orders of the sections of $f$ at $f^{-1}(P)$. Then $a_0=0$, and clearly
\[
a_i \geq m_i 
\]
for every $i=1, \dots, n$, so we see that $C$ ramifies at $P$ to order at least
\[
r_P= \sum_{i=1}^n (m_i-i).
\] 
As the codimension of Schubert varieties associated to the ramification at distinct points of $\mb{P}^1$ is additive, it suffices to prove that 
\begin{equation}\label{local_ramif_ineq}
r_P-1 \geq (n-2)g
\end{equation}
where the -1 on the left hand side arises from varying the preimage of $P$ along $\mb{P}^1$. 

We will show that \eqref{local_ramif_ineq} holds whenever $W_{\rm S} \leq 2g-1$.

\medskip
To this end, we use the fact \cite{BdM} that any numerical semigroup ${\rm S}$ of genus $g$ may be represented as a {\it Dyck path} $\tau=\tau({\rm S})$ on a $g \times g$ square grid $\La_{\rm S}$ with axes labeled by $0,1, \dots, g$. Each path starts at $(0,0)$, ends at $(g,g)$, and has unit steps upward or to the right. Namely, the $i$th step of $\tau$ is up if $i \notin {\rm S}$, and is to the right otherwise. The weight $W=W_{\rm S}$ of ${\rm S}$ is then equal to the total number of boxes in the Young tableau traced between the upper and left borders of the grid and the Dyck path $\tau$.

\medskip
Denote by $W^n_{\rm S}$ the contribution of the first $n$ elements $m_1, \dots, m_n$ of ${\rm S}$. Diagrammatically, this is precisely the area above the Dyck path in the subgrid determined by the first n columns of the $\La_{\rm S}$; equivalently,
\begin{equation}\label{W^n_S}
W^n_{\rm S}= ng- r_P.
\end{equation}
On the other hand, we clearly have 
\begin{equation}\label{W_S_versus_W^n_S}
W^n_{\rm S} \leq W_{\rm S} \leq 2g-1.
\end{equation}
Combining \eqref{W^n_S} with \eqref{W_S_versus_W^n_S} now yields \eqref{local_ramif_ineq}.



\subsection{The semigroup tree}
It is well-known that all numerical semigroups are indexed by the vertices of an infinite tree. More precisely, this {\it semigroup tree} has vertices indexed by sets of minimal generators;  the unique common ancestor of all vertices is the semigroup $\langle 1 \rangle= \mb{N}$ of genus 0. Any semigroup ${\rm S}$ of genus $g \geq 0$ may be presented by minimal generators as ${\rm S}=\langle a_1, \dots, a_N; b_1, \dots, b_M \rangle$ in which $a_1 < \dots < a_N < b_1 < \dots < b_M$ and the {\it conductor} $c=c({\rm S})$ satisfies $a_N<c \leq b_1$. (Here the set $\{b_1, \dots, b_M\}$ might be empty; recall that the conductor of a semigroup ${\rm S} \sub \mb{N}^r$ is the unique minimal element $c$ for which $c+ \mb{N}^r \sub {\rm S}$.) A child of ${\rm S}$, if it exists, is the unique completion of some subset ${\rm S}- \{b_i\}$, $i=1, \dots, M$ to a semigroup of genus $g+1$; for an illustration of the tree showing all vertices of genus $g \leq 4$, see \cite[Fig. 1]{BB}.

\medskip
The infinite branches of the semigroup tree distinguish infinite subfamilies of singularities, beginning with {\it hyperelliptic} singularities, i.e., those whose semigroups contain 2. Using the semigroup tree, it is easy (though tedious) to check that all {\it nonhyperelliptic} numerical semigroups ${\rm S}$ of genus $g<7$ verify $W_{\rm S} \leq 2g-1$.


\subsection{Conditions beyond ramification}
In this subsection, we consider {\it nonhyperelliptic} singularities whose semigroups ${\rm S}$ {\it fail} to verify the weight restriction $W_{\rm S} \leq 2g-1$. Specifically, we will show how to produce additional conditions in such cases, using multiplicative constraints imposed by ${\rm S}$. In this subsection we will do so for exactly those semigroups of genus $g \leq 8$ for which $W_{\rm S} \geq 2g$.

\begin{itemize}
\item {\bf Case: ${\rm S}= \langle 3,8 \rangle$}. 
The semigroup has genus 7, weight 14, and gap set $G_{\rm S}= \{1,2,4,5,7,10,13\}$. To estimate the codimension of the corresponding parameter space $\mc{V}_{\rm S}$ of morphisms $f$ with singularities of semigroup ${\rm S}$, we work analytically locally near such a singularity $P$. We may also assume that the lowest vanishing orders of global sections of $f$ in $f^{-1}(P)$ match the initial entries in ${\rm S}^*=\{0,3,6,8,9,11,12,14\}$. Indeed, it is straightforward to check that the space of morphisms that fail this minimality hypothesis is of codimension at least equal to $(n-2)g+1$. (In fact, the only cases in which ramification fails to produce the required $(n-2)g+1$ conditions are those in which $n \geq 6$.) It suffices to exhibit a single condition beyond (and independent of) ramification. Accordingly, let
\[
f_1(t)= t^3+ \al_1 t^4+ O(t^5) \text{ and } f_2(t)= t^6+ \al_2 t^7+ O(t^8)
\]
denote power series representatives for two local sections of minimal vanishing orders. Here we use local affine coordinates with respect to which $P$ is the origin in $\mb{C}^n$, and $t$ is a local coordinate centered in $f^{-1}(P)$. We then have
\[
f_1^2- f_2= (2 \al_1- \al_2) t^7+ O(t^8).
\]
In particular, since $7 \in G_{\rm S}$, it follows necessarily that $2 \al_1- \al_2=0$. So we obtain an additional linear condition on the coefficients of $f$ that is distinct from (and independent of) the ramification conditions.

\medskip

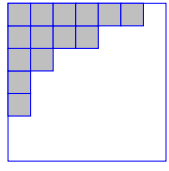
\begin{figure}

\begin{tikzpicture}[scale=0.30]
\draw[blue, very thin] (0,0) rectangle (7,7);
\filldraw[draw=blue, fill=lightgray] (0,2) rectangle (1,3);
\filldraw[draw=blue, fill=lightgray] (0,3) rectangle (1,4);

\filldraw[draw=blue, fill=lightgray] (0,4) rectangle (1,5);
\filldraw[draw=blue, fill=lightgray] (1,4) rectangle (2,5);

\filldraw[draw=blue, fill=lightgray] (0,5) rectangle (1,6);
\filldraw[draw=blue, fill=lightgray] (1,5) rectangle (2,6);
\filldraw[draw=blue, fill=lightgray] (2,5) rectangle (3,6);
\filldraw[draw=blue, fill=lightgray] (3,5) rectangle (4,6);

\filldraw[draw=blue, fill=lightgray] (0,6) rectangle (1,7);
\filldraw[draw=blue, fill=lightgray] (1,6) rectangle (2,7);
\filldraw[draw=blue, fill=lightgray] (2,6) rectangle (3,7);
\filldraw[draw=blue, fill=lightgray] (3,6) rectangle (4,7);
\filldraw[draw=blue, fill=lightgray] (4,6) rectangle (5,7);
\filldraw[draw=blue, fill=lightgray] (5,6) rectangle (6,7);

\end{tikzpicture}

\caption{Dyck path and Young tableau corresponding to ${\rm S}= \langle 3,8 \rangle$.}
\end{figure}

\item {\bf Case: ${\rm S}= \langle 3,10,17 \rangle$ }. Here $g({\rm S})=8$, $W_{\rm S}=16$, and $G_{\rm S}= \{1,2,4,5,7,8,11,14\}$. As in the preceding case, we work analytically in an adapted local coordinate $t$ and we may assume that local vanishing orders of $f$ in the preimage of the singularity are minimal. It suffices to exhibit a single condition beyond ramification. Once more, the linear  condition $2 \al_1- \al_2=0$ arises from the quadratic polynomial $f_1^2-f_2$ in the two local sections $f_1=t^3+ \al_1 t^4+ O(t^5)$ and $f_2=t^6+ \al_2 t^7+ O(t^8)$ of minimal vanishing orders.

\medskip
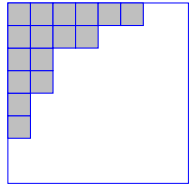
\begin{figure}

\begin{tikzpicture}[scale=0.30]
\draw[blue, very thin] (0,0) rectangle (8,8);
\filldraw[draw=blue, fill=lightgray] (0,2) rectangle (1,3);
\filldraw[draw=blue, fill=lightgray] (0,3) rectangle (1,4);

\filldraw[draw=blue, fill=lightgray] (0,4) rectangle (1,5);
\filldraw[draw=blue, fill=lightgray] (1,4) rectangle (2,5);

\filldraw[draw=blue, fill=lightgray] (0,5) rectangle (1,6);
\filldraw[draw=blue, fill=lightgray] (1,5) rectangle (2,6);

\filldraw[draw=blue, fill=lightgray] (0,6) rectangle (1,7);
\filldraw[draw=blue, fill=lightgray] (1,6) rectangle (2,7);
\filldraw[draw=blue, fill=lightgray] (2,6) rectangle (3,7);
\filldraw[draw=blue, fill=lightgray] (3,6) rectangle (4,7);

\filldraw[draw=blue, fill=lightgray] (0,7) rectangle (1,8);
\filldraw[draw=blue, fill=lightgray] (1,7) rectangle (2,8);
\filldraw[draw=blue, fill=lightgray] (2,7) rectangle (3,8);
\filldraw[draw=blue, fill=lightgray] (3,7) rectangle (4,8);
\filldraw[draw=blue, fill=lightgray] (4,7) rectangle (5,8);
\filldraw[draw=blue, fill=lightgray] (5,7) rectangle (6,8);

\end{tikzpicture}

\caption{Dyck path and Young tableau corresponding to ${\rm S}= \langle 3,10,17 \rangle$.}
\end{figure}

\item {\bf Case: ${\rm S}= \langle 4,6,13 \rangle$}. Here $g({\rm S})=8$, $W_{\rm S}=17$, and $G_{\rm S}= \{1,2,3,5,7,9,11,15\}$. For dimension reasons, we may assume local vanishing orders of sections in the preimage of the singularity are minimal in ${\rm S}^*=\{0,4,6,8,10,12,13,14,16\}$. We may further assume that $n \geq 6$. It will suffice to exhibit two linear conditions beyond ramification. For this purpose, consider three local (analytic) sections of minimal vanishing orders:
\[
f_1(t)= t^4+ \al_1 t^5+ O(t^6), f_2(t)= t^6+ \al_2 t^7+ O(t^8), \text{ and } f_3(t)= t^8+ \al_3 t^9+ O(t^{10}).
\]
We then have
\[
f_1^2-f_3= (2 \al_1-\al_3)t^9+ O(t^{10}) \text{ and } f_1^2f_2-f_2f_3= (2 \al_1-\al_2-\al_3)t^{15}+ O(t^{16}).
\]
Since 9 and 15 belong to $G_{\rm S}$, it follows that
\[
2 \al_1-\al_3=0 \text{ and } 2\al_1-\al_2-\al_3=0,
\]
which gives the two required conditions.

\begin{figure}

\begin{tikzpicture}[scale=0.30]
\draw[blue, very thin] (0,0) rectangle (8,8);
\filldraw[draw=blue, fill=lightgray] (0,3) rectangle (1,4);

\filldraw[draw=blue, fill=lightgray] (0,4) rectangle (1,5);
\filldraw[draw=blue, fill=lightgray] (1,4) rectangle (2,5);

\filldraw[draw=blue, fill=lightgray] (0,5) rectangle (1,6);
\filldraw[draw=blue, fill=lightgray] (1,5) rectangle (2,6);
\filldraw[draw=blue, fill=lightgray] (2,5) rectangle (3,6);

\filldraw[draw=blue, fill=lightgray] (0,6) rectangle (1,7);
\filldraw[draw=blue, fill=lightgray] (1,6) rectangle (2,7);
\filldraw[draw=blue, fill=lightgray] (2,6) rectangle (3,7);
\filldraw[draw=blue, fill=lightgray] (3,6) rectangle (4,7);

\filldraw[draw=blue, fill=lightgray] (0,7) rectangle (1,8);
\filldraw[draw=blue, fill=lightgray] (1,7) rectangle (2,8);
\filldraw[draw=blue, fill=lightgray] (2,7) rectangle (3,8);
\filldraw[draw=blue, fill=lightgray] (3,7) rectangle (4,8);
\filldraw[draw=blue, fill=lightgray] (4,7) rectangle (5,8);
\filldraw[draw=blue, fill=lightgray] (5,7) rectangle (6,8);
\filldraw[draw=blue, fill=lightgray] (6,7) rectangle (7,8);

\end{tikzpicture}

\caption{Dyck path and Young tableau corresponding to ${\rm S}= \langle 4,6,13 \rangle$.}
\end{figure}
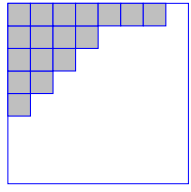

\end{itemize}

\begin{rem} In light of the preceding discussion, it is natural to ask for $g$-asymptotic estimates for the preponderance of semigroups of weight at least $2g$, i.e. refinements of Zhai's result \cite{Zh} that the number $n_g$ of {\it all} semigroups of given genus has Fibonacci-like asymptotics. Results of this type have recently been obtained by Kaplan and Ye \cite{KY}, who show among other things that the proportion of genus-$g$ semigroups with weight between (approximately) $.035g^2$ and $.089g^2$ is asymptotically equal to 1. In particular, this shows that semigroups of $g$-linearly bounded weight form a set of measure zero for $g \gg 0$.
\end{rem}

\subsection{Dimension counts for rational curves with hyperelliptic singularities}\label{hyp_dim_counts}
There is a unique hyperelliptic semigroup ${\rm S}^h= {\rm S}^h(g)$ of genus $g$, namely $\langle 2,2g+1 \rangle$, and it has weight $W_{{\rm S}^h}= \binom{g}{2}$. In particular, we have $W_{{\rm S}^h} \geq 2g$ for $g \geq 5$. As in the preceding section, we must use the arithmetic structure of ${\rm S}^h$ to produce conditions beyond ramification. It will also be convenient to codify the notion of the {\it ramification} arising from a numerical semigroup itself (as opposed to from one of its truncations). In general, given a numerical semigroup presented as 
\[
{\rm S}=\{s_i\}_{i \geq 0}
\]
where $s_i<s_{i+1}$ for all $i \geq 0$ and $s_0=0$, we denote by
\[
R_{\rm S}:= \{s_i-i\}_{i \geq 0} \text{ and } TR_{\rm S}:= \bigg\{\sum_{j=0}^i (s_j-j)\bigg\}
\]
the {\it ramification} and {\it total ramification} sequences of ${\rm S}$, respectively. We now specialize to the case in which ${\rm S}$ is hyperelliptic. It is then easy to see that $R_{{\rm S}^h}=\{r_i\}_{i \geq 0}$ and $TR_{{\rm S}^h}=\{\wt{r}_i\}_{i \geq 0}$ are characterized by
\[
\begin{split}
r_i&= i, 0 \leq i \leq g; \hspace{10pt} r_i=g, i \geq g; \text{ and }\\
\wt{r}_i&= \binom{i+1}{2}, 0 \leq i \leq g; \hspace{10pt} \wt{r}_i=\binom{g+1}{2}+ (i-g)g= ig- \binom{g}{2}, i \geq g.
\end{split}
\]
In particular, we see that
\begin{equation}\label{ramif_threshold}
\wt{r}_i \leq (i-2)g \hspace{10pt} \text{ if and only if } \hspace{10pt} g \geq \bigg\lfloor \frac{(i+1)i}{2(i-2)} \bigg\rfloor.
\end{equation}
On the other hand, the number of conditions beyond ramification required (to match $(n-2)g+1$ total conditions) is bounded above by $N_R=N_R(g)$, where
\begin{equation}\label{N_R}
N_R(g):=g(g-2)+1- \binom{g+1}{2}= \frac{g(g-5)}{2}+1.
\end{equation}

\begin{itemize}
\item {\bf Case: $g=5$.} Here \eqref{N_R} yields $N_R=1$, so we must produce a single condition beyond ramification. Moreover, applying \eqref{ramif_threshold}, we may assume $n \geq 4$. Accordingly, let
\[
f_1(t)= t^2+ \al_1 t^3+ O(t^4) \text{ and } f_2(t)= t^4+ \al_2 t^5+ O(t^6)
\]
denote two analytic sections of minimal nontrivial vanishing orders in a local coordinate $t$ adapted to the singularity. Then
\[
f_1^2-f_2= (2 \al_1- \al_2)t^5+ O(t^6).
\]
As $5 \in G_{{\rm S}^h}= \{1,3,5,7,9\}$, it follows that 
\begin{equation}\label{1sthypcondition}
2\al_1-\al_2=0
\end{equation}
which is the required condition.

\item {\bf Case: $g=6$.} Here $N_R=4$. More precisely, we must produce (at least) one, three, or four conditions beyond ramification when $n=3$, $n=4$, or $n \geq 5$, respectively. If $n=3$, we conclude using the same condition \eqref{1sthypcondition} as in the analysis of the $g=5$ case. If $n=4$, we further let
\[
f_3(t)= t^6+ \al_3 t^7+ O(t^8) \text{ and } f_4(t)= t^8+ \al_4 t^9+ O(t^{10}).
\]
Now consider the following quadratic polynomials in $f_1$, $f_2$, $f_3$, and $f_4$:
\[
f_1 f_3- f_4= (\al_1+ \al_3- \al_4)t^9+ O(t^{10}) \text{ and } f_2^2- f_4= (2\al_2-\al_4)t^9+ O(t^{10}).
\]
As $9 \in G_{{\rm S}^h}=\{1,3,5,7,9,11\}$, it follows that 
\begin{equation}\label{2dhypcondition}
\al_1+ \al_3- \al_4=0 \text{ and } 2 \al_2- \al_4=0.
\end{equation}
Taken together, conditions \eqref{1sthypcondition} and \eqref{2dhypcondition} allow us to conclude when $n=4$, provided that the global sections of the morphism $f$ that vanish in $f^{-1}(P)$ do so to minimal orders 2, 4, 6, and 8. The remaining possibility to be treated when $n=4$ is that the global sections of $f$ vanish to orders 2, 4, 6, and 10. However, in that instance only a single condition beyond ramification is required, and \eqref{1sthypcondition} is still operative.

\medskip
Now assume $n \geq 5$. We begin by considering the situation in which global sections of $f$ vanish to orders 0, 2, 4, 6, 8, and 10 in $f^{-1}(P)$. Accordingly, we let
\[
f_5(t)= t^{10}+ \al_5 t^{11}+ O(t^{12}).
\]
Then
\[
f_1 f_4- f_5= (\al_1+ \al_4- \al_5)t^{11}+ O(t^{12})
\]
and as $11 \in G_{{\rm S}^h}$, we deduce that
\begin{equation}\label{3dhypcondition}
\al_1+ \al_4- \al_5=0.
\end{equation}
The independent conditions \eqref{1sthypcondition}, \eqref{2dhypcondition}, and \eqref{3dhypcondition} allow us to conclude when $n \geq 5$, provided that the sections of $f$ vanish minimally to orders 0, 2, 4, 6, 8, and 10 in $f^{-1}(P)$. Similarly, if global sections of $f$ vanish minimally to orders 0, 2, 4, 6, 8, and $\ga$ with $\ga \geq 12$, we need only produce two conditions beyond ramification, and \eqref{1sthypcondition} and \eqref{2dhypcondition} remain operative. In all other cases, ramification gives the required $(n-2)g+1$ conditions. Note that conditions \eqref{1sthypcondition}, \eqref{2dhypcondition}, and \eqref{3dhypcondition} are all of the form
\begin{equation}\label{linearhypcondition}
\al_j= j \al_1
\end{equation}
where $j \geq 1$.

\item {\bf Case: $g=7$.} This time $N_R=8$. Specifically, we must produce two, five, seven, or eight conditions beyond ramification when $n=3$, $n=4$, $n=5$, or $n \geq 6$, respectively. If $n=3$, we may assume that sections of $f$ vanish minimally to orders 0, 2, 4, and 6 in $f^{-1}(P)$, since otherwise ramification yields $(n-2)g+1$ conditions. Imposing that the lowest-order terms of $f_1^2-f_2$ and $f_2^2-f_1f_3$ vanish gives the required two additional independent conditions, namely \eqref{1sthypcondition} and
\begin{equation}\label{4thhypcondition}
2\al_2- \al_1- \al_3=0.
\end{equation}
Of course, \eqref{1sthypcondition} and \eqref{4thhypcondition} simply translate to the instances $j=2,3$ of the linear constraint \eqref{linearhypcondition}.

\medskip
Now say $n=4$, and assume that sections of $f$ vanish minimally to orders 0, 2, 4, 6, and 8 in $f^{-1}(P)$. Much as before, we obtain three linear conditions \eqref{linearhypcondition}, $j=2,3,4$ beyond ramification. These may obtained by imposing that the leading terms of $F_1:=f_1^2-f_2$, $F_2:=f_1f_2-f_3$, and $F_3:=f_1f_3-f_4$, vanish. Thus $F_1$, $F_2$, and $F_3$, are power series with lowest terms of orders 6, 8, 10, and 12, respectively, in the absence of further nontrivial conditions. More precisely, rewriting
\[
f_j(t)= t^{2j}+ \al_j t^{2j+1}+ \be_j t^{2j+2}+ \ga_j t^{2j+3}+ O(t^{2j+4}), \hspace{10pt} j \geq 1
\]
to account for higher-order terms, we have
\[
\begin{split}
F_1&= (2\be_1+ \al_1^2- \be_2)t^6+ (2\ga_1+ 2\al_1 \be_1- \ga_2)t^7+ O(t^8), \\
F_2&= (\be_1+ \be_2+ \al_1 \al_2- \be_3)t^8+ (\ga_1+ \ga_2+ \al_1 \be_2+ \al_2 \be_1- \ga_3)t^9+ O(t^{10}), \text{ and }\\
F_3&= (\be_1+ \be_3+ \al_1 \al_3- \be_4)t^{10}+ (\ga_1+ \ga_3+ \al_1 \be_3+ \al_3 \be_1- \ga_4)t^{11}+ O(t^{12}).
\end{split}
\]
Imposing that
\[
[t^7](F_1- ([t^6]F_1)f_3)= [t^9](F_2- ([t^8]F_2)f_4)=0
\]
yields two further independent conditions as required, namely
\begin{equation}\label{1stnonlinhyp}
\begin{split}
2\ga_1+ 2\al_1\be_1- \ga_2- \al_3(2\be_1+ \al_1^2- \be_2) &= 0 \hspace{10pt} \text{ and } \\
\ga_1+ \ga_2+ \al_1 \be_2+ \al_2 \be_1- \ga_3- \al_4(\be_1+ \be_2+ \al_1 \al_2- \be_3) &=0.
\end{split}
\end{equation}
Note that the independence of these conditions is witnessed by the appearance of the variables $\ga_2=[t^7]f_2$ and $\ga_3=[t^7]f_3$ in the first and second of these equations, respectively, which in turn comes about because $f_2$ and $f_3$ are the ``linear" terms in their respective expansions
\[
F_1- ([t^6]F_1)f_3= f_1^2-f_2- ([t^6]F_1)f_3 \text{ and } F_2- ([t^8]F_2)f_4= f_1f_2-f_3- ([t^8]F_2)f_4.
\]
Finally, if instead sections of $f$ vanish to orders $(0,2,4,6,10)$, we need only exhibit three conditions beyond ramification, and it is easy to see that \eqref{linearhypcondition} holds for $j=2,3,5$. Similarly, if sections of $f$ vanish to orders $(0,2,4,8,10)$ or $(0,2,4,6,12)$, only one condition beyond ramification is required, and this may be chosen to be of type \eqref{linearhypcondition}. In all other cases, ramification furnishes the required number of conditions.

\medskip
The analysis is similar when $n=5$. Assume that sections of $f$ vanish to orders $(0,2,4,6,8,10)$. We must then produce seven conditions beyond ramification. It is easy to exhibit four of linear type, namely \eqref{linearhypcondition} with $2 \leq j \leq 5$. We obtain three further (independent, but nonlinear) conditions by imposing
\begin{equation}\label{nonlinhypbis}
[t^7](F_1- ([t^6]F_1)f_3)= [t^9](F_2- ([t^8]F_2)f_4)= [t^{11}](F_3-([t^{10}]F_3)f_5)=0.
\end{equation}
The argument is analogous, but simpler, if sections of $f$ vanish to orders larger than $(0,2,4,6,8,10)$.

\medskip
Finally, say $n \geq 6$. Assume that sections of $f$ vanish to orders $(0,2,4,6,8,12)$. We must produce eight conditions beyond ramification. Of these, five will be of linear type, namely \eqref{linearhypcondition} with $2 \leq j \leq 6$. The same three further nonlinear conditions \eqref{nonlinhypbis} as in the $n=5$ case remain operative.

\item {\bf Case: $g=8$.} Here $N_R=13$. More precisely, we must produce three, seven, ten, twelve, or thirteen conditions beyond ramification depending upon whether $n=3$, $n=4$, $n=5$, $n=6$, or $n \geq 7$.
We assume that sections of $f$ vanish to minimal possible orders, the remaining possibilities being analogous (and easier). Once the $(n-1)$ linear conditions \eqref{linearhypcondition} with $2 \leq j \leq n$ have been taken into account, we are left to produce one, four, six, seven, or seven nonlinear conditions. If $n=3$, we impose $[t^7](F_1- ([t^6]F_1)f_3)=0$, which we interpret as a condition on $[t^7]f_2$. Similarly, if $n=4$, we let
\[
Q_1:= F_1- ([t^6]F_1)f_3 \text{ and } Q_2:= F_2-([t^8]F_2)f_4.
\]
We impose
\begin{equation}\label{nonhypbis2}
[t^7]Q_1= [t^9]Q_2= 0;
\end{equation}
the first is the condition on $[t^7]f_2$ obtained in the $n=3$ case, while the second is a condition on $[t^9]f_3$. Now let
\[
\wt{Q}_1:= Q_1- ([t^8]Q_1)f_4 \text{ and } \wt{Q}_2:=Q_2-([t^{10}]Q_2)f_1f_4.
\]
We then impose
\begin{equation}\label{nonhypbis3}
[t^9]\wt{Q}_1= [t^{11}]\wt{Q}_2=0
\end{equation}
which give conditions on $[t^9]f_2$ and $[t^{11}]f_3$, respectively. Taken together, \eqref{nonhypbis2} and \eqref{nonhypbis3} give the required four nonlinear conditions.

\medskip
Now say $n=5$. Let $Q_3:= F_3- ([t^{12}]F_3)f_1f_5$; we impose that
\begin{equation}\label{nonhypbis4}
[t^{13}]Q_3=0;
\end{equation}
this is a condition on $[t^{13}]f_4$. Further, let $\wt{Q}_3:= Q_3- ([t^{14}]Q_3)f_3f_4$; we then impose that
\begin{equation}\label{nonhypbis5}
[t^{15}]\wt{Q}_3=0
\end{equation}
which is a condition on $[t^{15}]f_4$.
The conditions \eqref{nonhypbis2}, \eqref{nonhypbis3}, \eqref{nonhypbis4}, and \eqref{nonhypbis5} give the required six nonlinear conditions.   

\medskip
Finally, say $n \geq 6$. The six nonlinear conditions obtained in the $n=5$ case remain operative, and we obtain an additional condition, by setting $Q_4:= F_4- ([t^{14}]F_4)f_1f_6$ where $F_4:=f_1f_5-f_6$ (having implicitly imposed that $[t^{13}]F_4=0$, this being one of our linear conditions), and imposing that
\begin{equation}\label{nonhypbis6}
[t^{15}]Q_4=0
\end{equation}
which is a condition on $[t^{15}]f_6$.
\end{itemize}

\medskip
\begin{rems} We have not attempted here to make a systematic count of constraints on the coefficients of a morphism $f: \mb{P}^1 \ra \mb{P}^n$ arising from a hyperelliptic singularity in its image. In the follow-up paper \cite{CLM}, we show that a careful implementation of the basic strategy for obtaining conditions on coefficients $[t^{\rho}]f_j$ for $\rho \in G_{{\rm S}^h}$ as carried out when $g \leq 8$ above proves that morphisms with semigroup ${\rm S}^h= \langle 2, 2g+1 \rangle$ of arbitrary genus are of codimension {\it at least} $(n-1)g$; we further conjecture that their codimension is precisely the latter value. 

\medskip 
Note that when $n=g-1$, the number of conditions beyond ramification required to validate our original heuristic is exactly $N_R$. If $n=g-1$, our morphism $f: \mb{P}^1 \ra \mb{P}^n$ has a unique singularity which is hyperelliptic, and if the image of $f$ is not (globally) hyperelliptic, then $f$ describes a {\it canonical model} for $f$. This is a consequence of the Gorensteinness of a hyperelliptic singularity, which as explained in \cite{BDF} is manifested by the symmetry
\[
m_i= 2g-1- \ell_{g-i}
\]
between values $m_i$ and gaps $\ell_{g-i}$ of the semigroup ${\rm S}^h$ for every $i=0,\dots, g-1$, where by convention we set $m_0=0 \in {\rm S}^h$.
\end{rems}


\subsection{End of the proof of Thm~\ref{unibranch_theorem}}
The arguments of the preceding subsections show that Thm~\ref{unibranch_theorem} holds for all rational curves with a single cusp, so in this subsection we focus on rational curves with multiple cusps. The dimensional transversality result of Eisenbud--Harris is extremely useful in this regard, as it establishes that that ramification conditions in distinct points of (the domain of) a morphism $f: \mb{P}^1 \ra \mb{P}^n$ are independent. However, because of the presence of conditions beyond ramification, it isn't quite sufficient. Here we analyze the finite number of corresponding special cases, which we classify according to the partitions of $6 \leq g \leq 8$ induced by the the local genera of their cusps.

\begin{itemize}
\item {\bf Case: $g=6$.} Without loss of generality, we may assume the image of $f$ is of type $(5,1)$, where the entries of the partition refer to the genera of the singularities in the image of $f$. Further, we may assume the singularities are supported in $P_1=(0:\cdots:0:1)$ and $P_2=(1: 0\cdots:0)$ in $\mb{P}^n$ and 
are the images of $(0:1)$ and $(1:0)$ in $\mb{P}^1$, respectively; and that the corresponding semigroups are $\langle 2,5 \rangle$ and $\langle 2,3 \rangle$, respectively. The independence of the (linear) conditions beyond ramification of the genus-5 hyperelliptic singularity in $P_1$ with respect to the simple cusp in $P_2$ is clear. Indeed, considering the matrix of coefficients of our degree-$d$ morphism as a $(n+1) \times (d+1)$ box, the vanishing conditions imposed by the singularities in $P_1$ and $P_2$ define diagonally opposite Young tableaux inside the box.
\item {\bf Case: $g=7$.} Without loss of generality, we may assume $f$ to be of partition type $(5,2)$, $(6,1)$ or $(5,1,1)$. Of these possibilities, the first two may be handled just as in the $g=6$ case. Similarly, if $f$ is of type $(5,1,1)$, we may assume its singularities are the images of $(0:1)$, $(1:0)$, and $(1:1)$, respectively. Likewise, we may assume that $f$ maps $(0:1)$ to $(0:\cdots:0:1)$ and $(1:0)$ to $(1: 0:\cdots:0)$. The remaining singularity, a simple cusp, is specified by requiring $\mbox{rank}(Df(1:1)) \leq 1$, where $Df=(\frac{\partial{f}}{\partial t},\frac{\partial{f}}{\partial u})$. 
To count conditions explicitly, we assume the orders of vanishing of $f$ in $(0:1)$, the preimage of the singularity of genus 5, are minimal (as usual, this is the crucial case). Applying a diagonal automorphism if necessary, we may then assume that $f=(f_0,\dots,f_n)$ is of the form
\[
\begin{split}
f_0 &= u^d+ a_{0,1} tu^{d-1}+ \dots+ a_{0,d-n-1} t^{d-n-1}u^{n+1} \\
f_1 &= t^2 u^{d-2}+ a_{1,3} t^3 u^{d-3}+ \dots+ a_{1,d-n} t^{d-n} u^n \\
&\dots \\
f_{n-1} &= t^{2n-2} u^{d-2n+2}+ a_{n-1,2n-1} t^{2n-1} u^{d-2n+1}+ \dots+ a_{n-1,d-2} t^{d-2} u^2 \\
f_n &= t^{2n} u^{d-2n}+ a_{n,2n+1} t^{2n+1} u^{d-2n-1}+ \dots+ a_{n,d} t^d
\end{split}
\]
with $a_{i,j} \in \mb{C}$ for all $0 \leq i \leq n$ and $0 \leq j \leq d$. Further, \eqref{1sthypcondition} yields $2(a_{1,3}-a_{0,1})=a_{2,5}-a_{0,1}$, i.e. 
\begin{equation}\label{conditionbeyondramif}
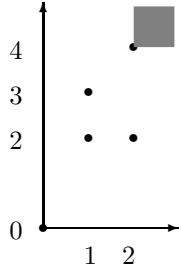
a_{0,1}= 2a_{1,3}-a_{2,5}.
\end{equation}
Note that
\begin{equation}\label{partial_t}
\frac{\partial f_j}{\partial t}(1,1)= \sum_{i=j}^{d-n-1} (i+j)a_{j,i+j} \text{ for } 0 \leq j \leq n-1, \text{ and } \frac{\partial f_n}{\partial t}(1,1)= \sum_{i=n}^{d-n-1} (n+i) a_{n,n+i}
\end{equation}
and
%
\begin{equation}\label{partial_u}
\frac{\partial f_j}{\partial u}(1,1)= \sum_{i=j}^{d-n-1} (d-i-j) a_{j,i+j} \text{ for } 0 \leq j \leq n-1, \text{ and } \frac{\partial f_n}{\partial u}(x,1)= \sum_{i=n}^{d-n-1} (d-n-i) a_{n,n+i}
\end{equation}
where $a_{0,0}=a_{1,2}=\dots=a_{n-1,2n-2}=a_{n,2n}:=1$. The requirement that $\mbox{rank}(Df(1:1)) \leq 1$ translates to
\begin{equation}\label{vanishingminors}
\frac{\partial f_j}{\partial t}(1,1) \frac{\partial f_k}{\partial u}(1,1)- \frac{\partial f_k}{\partial t}(1,1) \frac{\partial f_j}{\partial u}(1,1)=0
\end{equation}
for all $j \neq k \in \{0,\dots,n\}$. Using \eqref{partial_t} and \eqref{partial_u}, is clear that we obtain (in fact more than) the required $(n-2)$ independent conditions on the coefficients $a_{i,j}$ from the determinantal equations \eqref{vanishingminors}.
\item {\bf Case: $g=8$.} 
Invoking the ``box principle" of the preceding two items, we may assume without loss of generality that $f$ is associated to a partition with at least four nontrivial parts. The unique such possibility is 
$(5,1,1,1)$. But because each simple cusp imposes $(n-1)$ instead of the required $(n-2)$ ramification conditions (and because the genus-5 singularity is at worst hyperelliptic, in which case its ramification fails by at most a single condition to meet the codimension heuristic in its genus) it follows by Eisenbud--Harris dimensional transversality that ramification arising from singularities of type $(5,1,1,1)$ produce (more than) the required $8(n-2)$ conditions. (Note that this argument also serves in the $(5,1,1)$ case.)
\end{itemize}
\qed

\medskip
\begin{rem} Our study complements the largely topological investigations of (the existence of) {\it planar} rational cuspidal curves carried out, e.g., in \cite{deB}, \cite{FZ}, \cite{O}, and \cite{Ton}; see \cite{Mo} for an exhaustive bibliography. Piontkowski conjectured in \cite{Pi} that a rational plane curve of degree at least six may have at most three cusps, and Koras--Palka \cite{KoP} have announced that a proof of this result is forthcoming. It is natural to wonder whether a similar boundedness result holds for rational cuspidal curves in $\mb{P}^n$, $n \geq 3$. 
\end{rem}

\section{Severi-type varieties associated with $\ga^*$-hyperelliptic cusps}
Given a choice of nonnegative integers $d$, $g$, and $n$, we will refer to the space of morphisms $M^n_{d,g}$ as a {\it Severi-type} variety. When $n=2$, Harris' celebrated result \cite{H1} implies that $M^n_{d,g}$ is irreducible, and equal to the closure (inside $M^n_d$) of the subvariety of morphisms whose images contain $g$ nodes. It is natural to ask whether this phenomenon persists when $n \geq 3$. In this section, we will show that when $n \geq 8$ this is not the case, as $M^n_{d,g}$ contains ``excess" components of dimension strictly larger than that of the $g$-nodal sublocus.

\medskip
The excess components of Severi-type varieties that we construct arise as parameterizations $f(t)=(f_0(t),\dots,f_n(t))$ in $t$-power series whose associated $t$-th-power valuations $(\text{val}_t(f_0),\dots,\text{val}_t(f_n))$ arise from arithmetic progressions of a particular type, and whose higher-order terms have generic coefficients. The value semigroups of the cusps determined by these parameterizations are described by the following definition.

\begin{dfn}
Given a nonnegative integer $\ga$, we say that a numerical semigroup ${\rm S} \sub \mb{N}$ is $\ga^{\ast}$-hyperelliptic provided that
\begin{enumerate}
    \item $\ga$ even numbers belong to ${\rm S} \cap \{1,\dots,4\ga-2\}$; and
    \item $4\ga \in {\rm S}$.
\end{enumerate}
\end{dfn}
We call such semigroups {\it $\ga^{\ast}$-hyperelliptic} by virtue of their similarity to the $\ga$-hyperelliptic semigroups studied by Torres in \cite{To2}.

\subsection{$\ga^{\ast}$-hyperelliptic cusps}
Given a nonnegative integer $\ga$ as above, we now define $t$-power series 
\begin{equation}\label{parameterizing_functions}
f_i(t):=t^{2(\ga+i)}+O(t^{2(\ga+i)+1})
\end{equation}
for every $i=0,\dots,\ga-1$, in which that the coefficients of $t$-powers strictly larger that $2i$ are generically chosen. Then $f(t)=(f_0(t),\dots,f_{\ga-1}(t))$ parameterizes a unibranch singularity in $\mb{C}^{\ga}$ whose associated value semigroup ${\rm S^{\ast}}$ is $\ga^{\ast}$-hyperelliptic.

\begin{prop}\label{counterexample_semigroup}
The value semigroup ${\rm S^{\ast}}$ determined by the parameterizing functions \eqref{parameterizing_functions} decomposes as ${\rm S}^{\ast}={\rm S}^{\ast}_0+ {\rm S}^{\ast}_1$, where 
\[
{\rm S}^{\ast}_0=\langle 2\ga, 2\ga+2, \dots, 4\ga-2 \rangle \text{ and } {\rm S}^{\ast}_1= \langle 4\ga+5, 4\ga+7, \dots, 6\ga+3 \rangle.
\]
\end{prop}

\begin{proof}
Clearly ${\rm S}^{\ast}_0+ {\rm S}^{\ast}_1 \sub {\rm S}^{\ast}$. Indeed, ${\rm S}^{\ast}_0$ is generated by the $t$-valuations of the parameterizing functions $f_i$, while ${\rm S}^{\ast}_1$ is generated by $t$-valuations of {\it quadratic} binomials in the $f_i$. Moreover, for each $i=0,1$, it is clear that every element in ${\rm S}^{\ast}_i$, has $2$-residue equal to $i$, and that ${\rm S}^{\ast}_0$ (resp., ${\rm S}^{\ast}_1$) contains every even (resp., odd) element in $\mb{N}$ greater than or equal to $2\ga$ (resp., $4\ga+5$). We conclude by observing that $2\ga$ (resp., $4\ga+5$) is the {\it minimal} positive even (resp., odd) integer that belongs to ${\rm S}^{\ast}$.
\end{proof}

The upshot of Proposition~\ref{counterexample_semigroup} is that the gap set $G_{{\rm S}^{\ast}}$ associated with ${\rm S^{\ast}}$ decomposes as $G_{{\rm S}^{\ast}}= G_0 \sqcup G_1$, where
\[
G_0=\{2,4,\dots,2\ga-2\} \text{ and } G_1=\{1,3, \dots, 4\ga+3\}.
\]
In particular, the value semigroup ${\rm S}^{\ast}$ is of genus 
\begin{equation}\label{genus_calc}
g({\rm S}^{\ast})=\#G_0+ \#G_1=3\ga+1.
\end{equation}
On the other hand, by construction the only conditions imposed on rational curves by singularities parameterized by \eqref{parameterizing_functions} arise from ramification, and there are
\begin{equation}\label{ramif_calc}
r=r(f)= ((2\ga-1)+(4\ga-2-\ga)) \cdot \ga/2-1= \frac{(5\ga-3)\ga}{2}-1
\end{equation}
of these. Note here that we have subtracted one in \eqref{ramif_calc} to account for variation in the preimage of the cusp parameterized by \eqref{parameterizing_functions}. It follows immediately from equations \eqref{genus_calc} and \eqref{ramif_calc} that
\[
r(f) < (\ga-2)g({\rm S}^{\ast})
\]
for all $\ga \geq 8$. Our discussion proves the following result.

\begin{thm}\label{severi_counterexamples}
For every positive integer $n \geq 8$, there exist Severi-type varieties $M^n_{d,g}$ with $g=3n+1$ and $d \gg g$ that contain excess components not contained in the closure of the space of $g$-nodal rational curves of degree $d$ in $\mb{P}^n$. In particular, these Severi-type varieties are reducible.
\end{thm}

Theorem~\ref{severi_counterexamples} leads naturally to the question of whether such excess components admit interesting {\it global} geometric interpretations; for example, at the level of the equations of their images.  
It would also be interesting to know whether generalizations of the $\ga^{\ast}$-hyperelliptic semigroups introduced above yield more general constructions of reducible Severi-type varieties; in particular, it is natural to ask whether the hypothesis $n \geq 8$ in Theorem~\ref{severi_counterexamples} may be slackened to $n \geq 3$.


\begin{thebibliography}{30}
\bibitem{ACGH} E. Arbarello, M. Cornalba, P. A. Griffiths, and
  J. Harris, ``Geometry of Algebraic Curves'', Springer, 1985.
\bibitem{BDF} V. Barucci, M. D'Anna, and R. Fr\"oberg, {\it Analytically unramified one-dimensional semi local rings and their value semigroups}, J. Pure Appl. Alg. {\bf 147} (2000), 215--254.
\bibitem{BF} V. Barucci, and R. Fr\"oberg, {\it One-dimensional almost Gorenstein rings}, J. Algebra {\bf 188} (1997), no. 2, 418--442.
\bibitem{BdM} M. Bras-Amor\'os and A. de Mier, {\it Representation of numerical semigroups by Dyck paths}, Semigroup Forum {\bf 75} (2007), no. 3, 676--681.
\bibitem{BB} M. Bras-Amor\'os and S. Bulygin, {\it Towards a better understanding of the semigroup tree}, Semigroup Forum {\bf 79} (2009), no. 3, 561--574.
\bibitem{BIV} J. Buczy\'nski, N. Ilten, and E. Ventura, {\it Singular curves of low degree and multifiltrations from osculating spaces}, \url{arXiv:1905.11860}.
\bibitem{CH} E. Carvalho and M.E. Hernandes, {\it The semiring of values of an algebroid curve}, \url{arXiv:1704.04948}.
\bibitem{CS} A. Contiero and K.-O. St\"ohr, {\it Upper bounds for the dimension of moduli spaces of curves with symmetric Weierstrass semigroups}, J. London Math. Soc. {\bf 88} (2013), 580--598.
\bibitem{Co1} E. Cotterill, {\it Rational curves of degree 11 on a general quintic threefold}, Quart. J. Math. {\bf 63} (2012), no. 3, 539--568.
\bibitem{Co} E. Cotterill, {\it Rational curves of degree 16 on a general heptic fourfold}, J. Pure Appl. Alg. {\bf 218} (2014), 121--129.
\bibitem{CFM} E. Cotterill, L. Feital, and R.V. Martins, {\it Singular rational curves in $\mb{P}^n$ via semigroups}, \url{arXiv:arXiv:1511.08515}.
\bibitem{CLM} E. Cotterill, V. Lima, and R.V. Martins, {\it Rational curves with $\ga$-hyperelliptic singularities}, in preparation.
\bibitem{CFHM} E. Cotterill, L. Feital, M.E. Hernandes, and R.V. Martins, {\it Dimension counts for multibranch singular rational curves in $\mb{P}^n$ via semigroups}, in preparation.
\bibitem{CKPU} D. Cox, A. Kustin, C. Polini, and B. Ulrich, {\it A study of singularities on rational curves via syzygies}, Mem. Amer. Math. Soc. {\bf 222} (2013), no. 1045.
\bibitem{CT} C. Carvalho and F. Torres, {\it On numerical semigroups associated to coverings of curves}, Semigroup Forum {\bf 67} (2003), no. 3, 344--354.
\bibitem{EH} D. Eisenbud and J. Harris, {\it Divisors on general curves and cuspidal rational curves}, Invent. Math. {\bf 74} (1983), 371--418.
\bibitem{deB} J. Fern\'andez de Bobadilla, I. Luengo, A. Melle-Hern\'andez, and A. N\'emethi, {\it On rational cuspidal plane curves, open surfaces, and local singularities}, pp.411--442 in ``Singularity theory", World Sci. Publ., 2007.
\bibitem{FZ} H. Flenner and M. Zaidenberg, {\it Rational cuspidal plane curves of type $(d,d-3)$}, Math. Nachr. {\bf 210} (2000), 93--110.
\bibitem{GH} P. A. Griffiths and J. Harris, {\it On the variety of special linear systems on a general algebraic curve}, Duke Math. J. {\bf 47} (1980), no. 1, 233--272.
\bibitem{H1} J. Harris, {\it On the Severi problem}, Invent. Math. {\bf 84} (1986), vol. 3, 445--461. 
\bibitem{HM} J. Harris and I. Morrison, ``Moduli of curves", Springer, 1998.
\bibitem{JK} T. Johnsen and S. Kleiman, {\it Rational curves of degree at most 9 on a general quintic threefold}, Comm. Alg. {\bf 24} (1996), no. 8, 2721--2753.
\bibitem{KY} N. Kaplan and L. Ye, {\it The proportion of numerical semigroups}, J. Alg. {\bf 373} (2013), 377--391.
\bibitem{KoP} M. Koras and K. Palka, {\it The Coolidge--Nagata conjecture}, \url{arXiv:1502.07149}.
\bibitem{Lau} G. Laumon, {\it Fibres de Springer et jacobiennes compactifi\'ees}, in ``Algebraic geometry and number theory," Prog. Math. {\bf 253} (2006), 515--563.
\bibitem{Mo} T. Karoline Moe, {\it Rational cuspidal curves}, Univ. Oslo master's thesis, 2008; \url{arXiv:1511.02691}.
\bibitem{O} S. Y. Orevkov, {\it On rational cuspidal curves I. Sharp estimate for degree via multiplicities}, Math. Ann. {\bf 324} (2002), no. 4, 657--673.
\bibitem{Pi} J. Piontkowski, {\it On the number of cusps of rational cuspidal plane curves}, Exp. Math. {\bf 16} (2007), no. 2, 251--255.
\bibitem{Ton} K. Tono, {\it On the number of cusps of cuspidal plane curves}, Math. Nachr. {\bf 278} (2005), no. 1-2, 216--221.
\bibitem{To2} F. Torres, {\it On $\ga$-hyperelliptic numerical semigroups}, Semigroup Forum {\bf 55} (1997), 364--379.
\bibitem{Zh} A. Zhai, {\it Fibonacci-like growth of the number of semigroups of given genus}, Semigroup Forum {\bf 86} (2013), no. 3, 634--662.
\end{thebibliography}
\end{document}